\documentclass[letterpaper, 10 pt, conference]{ieeeconf}  

\IEEEoverridecommandlockouts                              

\usepackage{graphicx}

\usepackage{url}
\usepackage{amsmath,amssymb,mathtools}
 \usepackage{amsthm,thmtools}
\usepackage{amsfonts} 

\usepackage{enumerate}
\usepackage{color}                  

\usepackage{verbatim}

\usepackage{amssymb}
\usepackage{amsbsy}
\usepackage{amsmath}
\usepackage{setspace}
\usepackage{url}

\usepackage{subcaption}
\usepackage{float}
 \usepackage{dsfont}
 
\usepackage{tikz}
\usetikzlibrary{positioning}

 \newcommand{\diag}{\operatorname{diag}}

\newcommand{\spanop}{\operatorname{span}}

\newcommand*\diff{\mathop{}\!\mathrm{d}}

\newcommand{\trace}{\operatorname{trace}}

 \declaretheorem[qed=$\blacksquare$ ]{Example}
 \declaretheorem[name={Proof},    numbered=no,qed=$\blacksquare$ ]{Proof}

\declaretheorem[name={Definition}  ] {Definition}

\declaretheorem[name={Remark}  ] {Remark}
\declaretheorem[name={Corollary}  ] {Corollary}

\declaretheorem[name={Proposition}  ] {Proposition}

\newcommand {\R}{\mathbb R}

\newcommand{\be}{\begin{equation}}
\newcommand{\ee}{\end{equation}}

\usepackage{lineno}

\DeclareMathOperator{\rank}{rank}

 \title{\LARGE \bf
Serial interconnections of~$1$-contracting and~$2$-contracting~systems }

\author{Ron Ofir, Michael Margaliot, Yoash Levron, and
 Jean-Jacques Slotine
\thanks{This research was  partially supported by a research grant from the Israel Science Foundation~(ISF).}
\thanks{RO and YL are  with 
		the Andrew and Erna Viterbi Faculty of Electrical Eng., Technion---Israel Institute of Technology, Haifa 3200003, Israel.}%
\thanks{MM (Corresponding Author) is  with the School of Electrical  Eng.,
		and the Sagol School of Neuroscience, 
		Tel-Aviv University, Tel-Aviv~69978, Israel.
		E-mail: \texttt{michaelm@tauex.tau.ac.il}}%
		\thanks{JJS is with the Department of Mechanical Eng. and
        the Department of Brain and Cognitive Sciences, Massachusetts Institute of Technology, Cambridge, Massachusetts, USA.}
}

\begin{document}

\maketitle
\begin{abstract}
    The flow  of contracting   systems  contracts 1-dimensional polygons (i.e. lines) at an exponential rate.  One reason for the  usefulness of contracting  systems  is that many interconnections of contracting sub-systems yield an overall  contracting system. 
    A recent generalization of contracting systems  is called~$k$-contracting systems, where~$k\in\{1,\dots,n\}$. 
    The flow  of such systems contracts $k$-dimensional polygons at an exponential rate, and in particular they reduce to contracting systems when~$k=1$. 
    Here, we analyze serial interconnections of~$1$-contracting and~$2$-contracting systems. We provide conditions guaranteeing that such interconnections  have a well-ordered asymptotic behaviour, and demonstrate  the theoretical results  using several   examples.
\end{abstract}
\section{Introduction}
Contracting  systems have found 
numerous applications in systems and control theory. This is due to several reasons. 
First, contracting systems have a well-ordered behaviour: any two trajectories approach one another at an exponential rate~\cite{LOHMILLER1998683}.  
In particular, if an equilibrium point exists then it is unique and globally exponentially stable. If the vector field is~$T$-periodic then the system entrains, i.e. all solutions converge exponentially to a unique  $T$-periodic trajectory~\cite{LOHMILLER1998683,entrain2011}. In fact, contracting
systems have a well-defined frequency response, as shown in~\cite{freq_convergent} in the context of convergent  systems~\cite{pavlov_book}. 
Second, there exist simple sufficient conditions for contraction based on matrix measures~\cite{LOHMILLER1998683,sontag_contraction_tutorial}. Third, various interconnections of contracting systems, including parallel, serial, and feedback connections,  yield an overall contracting system~\cite{LOHMILLER1998683,netwok_contractive}. 

Ref.~\cite{kordercont} studied a generalization called~$k$-contraction   (see also the note~\cite{weak_manchester}), with~$k\in\{1,\dots,n\}$. The flow of such systems 
contracts~$k$-dimensional polygons at an exponential rate.  In particular, for~$k=1$ these are just standard contracting  systems. 
This generalization is motivated in part by the seminal work by
  Muldowney and his colleagues~\cite{muldo1990,li1995}, on systems that, using the new terminology, are~$2$-contracting in a constant metric. Roughly speaking, every bounded solution of a time-invariant  $2$-contracting system converges to an equilibrium point. This is different from the case of $1$-contracting systems,
  as the equilibrium point is not necessarily unique.

Contraction theory is an active area of research. 
Recent contributions include
contraction on  Riemannian manifolds~\cite{SIMPSONPORCO201474}, various notions of ``weak contraction'' (see, e.g.~\cite{jafarpour2020weak,cast_book}), contraction of piecewise-smooth  dynamical systems~\cite{contraction_non_diff},   analysis of learning algorithms using contraction theory~\cite{wensing}, and the introduction of~$\alpha$-contracting systems, with~$\alpha\geq 1$ \emph{real}, which is motivated in part by the seminal
works of Douady and   Oesterl\'{e}~\cite{Douady1980}, and 
Leonov  and  his   colleagues (see the recent monograph by Kuznetsov and   Reitmann~\cite{book_volker2021}) on bounding the Hausdorff dimension of complex attractors.

Since many interconnections of contracting systems yield an overall contracting system, it is natural to ask if the same holds for~$k$-contracting systems as well~\cite{weak_manchester}.
Here, we address this question in some detail for~$k$-contracting systems with~$k\in\{1,2\}$.
This problem is more delicate than in the case of $1$-contracting systems because the well-ordered behaviour of $2$-contracting systems
only holds in the time-invariant case, while connecting two systems implies that at least one system has an input from the other system  and thus is time-varying.

Our main contribution is a proof that various serial connections of $k$-contracting systems, with~$k \in \{1,2\}$, have a ``well-ordered'' asymptotic behaviour: they have no non-trivial periodic solutions, and, under stronger assumptions, all solutions converge to an equilibrium point (which is not necessarily unique). We also show that such connections are in general neither $1$-contracting nor $2$-contracting, and thus our results may be used to analyze systems that cannot be studied using only the theory of~$2$-contracting systems. To apply our results to wider set of systems, we also provide sufficient conditions guaranteeing that a given system can be decomposed as the serial connection of two systems.

The next section reviews known definitions and results that are used later on. Section~\ref{sec:main}
includes the main results, 
and the final section concludes. Due to space limitations, we focus on theoretical results and provide only a few applications in Section~\ref{sec:app}. More applications will appear in an extended version of this note 
that is now in preparation. 

We use standard notation. Small [capital]  letters denote column vectors [matrices]. $I_n$ is the~$n\times n$ identity matrix. 
For a matrix~$A$,~$A^T$ is the transpose of~$A$.
If~$A$ is square, then~$|A|$
[$\trace(A)$] is the determinant [trace] of~$A$.

\section{Preliminaries} 
The sufficient  condition for~$k$-contraction
    in~\cite{kordercont}
is based on the~$k$th additive compound of the Jacobian of the vector field. 
To make this note more accessible, we briefly review these topics. 
For more details, see also~\cite{muldo1990}. For more  recent applications of these compounds in systems and control theory, see~\cite{margaliot2019revisiting,rami_osci,pines2021,cheng_diag_stab,barshalom2021compound}.

Let~$C\in\R^{n\times m}$. For~$k\in\{1,\dots,\min\{n,m\}\}$, 
the   \emph{$k$th multiplicative compound} of~$C$, denoted~$C^{(k)} $, is the~$\binom{n}{k}\times\binom{m}{k}$ matrix that contains all the~$k\times k$ minors of~$C$ in lexicographic order~\cite{muldo1990}. For example, for~$n=m=3$ and~$k=2$,
$C^{(2)}$ is the~$3\times 3$ matrix:
\[
 \begin{bmatrix}
    \begin{vsmallmatrix}  c_{11}& c_{12} \\ c_{21}& c_{22} \end{vsmallmatrix}  &
         \begin{vsmallmatrix}  c_{11}& c_{13} \\ c_{21}& c_{23} \end{vsmallmatrix}  &
             \begin{vsmallmatrix}  c_{12}& c_{13} \\ c_{22}& c_{23} \end{vsmallmatrix}  \\
     \begin{vsmallmatrix}  c_{11}& c_{12} \\ c_{31}& c_{32} \end{vsmallmatrix}  &
         \begin{vsmallmatrix}  c_{11}& c_{13} \\ c_{ 31}& c_{33} \end{vsmallmatrix}  &
             \begin{vsmallmatrix}  c_{12}& c_{13} \\ c_{32}& c_{33} \end{vsmallmatrix}  \\
   \begin{vsmallmatrix}  c_{21}& c_{22} \\ c_{31}& c_{32} \end{vsmallmatrix}  &
         \begin{vsmallmatrix}  c_{21}& c_{23} \\ c_{31}& c_{33} \end{vsmallmatrix}  &
             \begin{vsmallmatrix}  c_{22}& c_{23} \\ c_{32}& c_{33} \end{vsmallmatrix}  
\end{bmatrix},
\]
where~$|B|$ denotes the determinant of~$B$.
In particular, $C^{(1)}=C$ and if~$n=m$ then~$C^{(n)}=|C|$.
The Cauchy–Binet  formula~\cite[Chapter~0]{matrx_ana},  asserts  that for any~$B \in\R^{n\times m},C\in\R^{m\times p}$ and any~$k\in \{1,\dots, \min \{n,m,p\}\}$,
we have
\be\label{eq:mcp}
(BC)^{(k)}=B^{(k)}  C^{(k )}.
\ee
This justifies the term multiplicative compound. 
In particular,~\eqref{eq:mcp} 
implies  that if~$n=m=p$ then~$|BC|=|B||C|$, and that if~$A\in \R^{n\times n}$ is   non-singular then~$(A^{(k)})^{-1}=(A^{-1})^{(k)}$.

Let~$A\in\R^{n\times n}$ with 
eigenvalues~$\lambda_i$, $i\in\{1,\dots,n\}$. The  eigenvalues of~$A^{(k)}$  are
$
   \lambda_ {i_1} \lambda_{i_2} \dots \lambda_{i_k  } ,
$
with~$1\leq i_1<i_2 <\dots <i_k\leq n   $.

For~$k\in\{1,\dots,n\}$, 
the   \emph{$k$th additive compound} of~$A$    is the~$\binom{n}{k}\times\binom{n}{k}$ matrix defined by
\[
A^{[k]}:=\frac{d}{d\varepsilon} (I+\varepsilon
 A)^{(k)} |_{\varepsilon=0}.
\]
In other words,~$(I+\varepsilon
 A)^{(k)}=I+\varepsilon A^{[k]} +o(\varepsilon)$. In particular,
 $
 A^{[1]}=A $, 
  and  $A^{[n]}=\trace(A).
 $
The eigenvalues of~$A^{[k]}$  are
$
   \lambda_ {i_1}+ \lambda_{i_2} +\dots+ \lambda_{i_k  } ,
$
with~$1\leq i_1<i_2 <\dots <i_k\leq n   $.

 It is useful to know how these compounds are affected by a coordinate transformation. 
 Let~$V\in\R^{m \times n}, W \in \R^{n\times p}$. Then~\eqref{eq:mcp}
 yields
\[
(V A W)^{(k)}=V^{(k)} A^{(k)} W^{(k)}. 
\]
If, in addition,~$m=p$ and~$V W = I_m$ then
\begin{align}\label{eq:ctf}
    (V A W)^{[k]} &= \frac{d}{ d \varepsilon} \left( I_m + \varepsilon  V A W \right)^{(k)}|_{\varepsilon=0} \nonumber \\
    &=  \frac{d}{d \varepsilon} ( V ( I_n + \varepsilon A)W )^{(k)}|_{\varepsilon=0} \nonumber \\
    &= V^{(k)} A^{[k]} W^{(k)}.
\end{align}

In the context of dynamical systems, the importance of these compounds  is due to following fact. If~$\Phi:\R _+\to \R^{n\times n}$ is the solution of the matrix differential equation
\[
\frac{d}{dt} \Phi(t)=A(t) \Phi(t),\quad \Phi(0)=I, 
\]
where~$t\to A(t)$ is continuous, then
\be\label{eq:rs}
\frac{d}{dt} (\Phi(t))^{(k)}=(A(t))^{[k]}
(\Phi(t))^{(k)}. 
\ee
In other words,~$  (\Phi(t))^{(k)}$
also evolves according to  a linear dynamics with the matrix~$(A(t))^{[k]}$.
Roughly speaking,~$\Phi^{(k)}$ determines the evolution of $k$-dimensional polygons  under the dynamics~$\dot x=Ax$~\cite{wini2010}.

Recall that a vector norm~$|\cdot|:\R^n\to\R_+$ induces a matrix norm~$||A||:=\max_{|x|=1} |Ax|$, and a matrix measure~$\mu(A):=\lim_{\varepsilon\to 0^+}  (||I+\varepsilon A||-1)/ \varepsilon$. If~$\mu( (A(t))^{[k]} )\leq-\eta<0$ all~$t\geq 0$ then
applying Coppel's inequality~\cite{coppel1965stability} to~\eqref{eq:rs} yields
 $|| (\Phi(t))^{(k)} ||\leq \exp(-\eta t) 
|| (\Phi(0))^{(k)} ||$ for all~$t\geq 0$. This leads to the following. 
\begin{Definition}\cite{kordercont}
Consider the nonlinear system~$\dot x(t)=f(t,x(t))$, with~$f$ a~$C^1$ mapping,
and suppose that its trajectories evolve on a convex set~$\Omega\subseteq \R^n$. Let~$J(t,x):=\frac{\partial}{\partial x} f(t,x)$ denote the Jacobian of~$f$ with respect to~$x$. The system is called~\emph{$k$-contracting} if
\be\label{eq:gpo}
\mu \left( (J(t,z))^{[k]} \right )\leq -\eta<0,\text{ for all }
t\geq 0,z\in \Omega.
\ee
\end{Definition}
Note that for~$k=1$ this reduces to the  standard infinitesimal contraction  condition~\cite{sontag_contraction_tutorial},
as~$J^{[1]}=J$. Note also that condition~\eqref{eq:gpo} is robust in the sense that if it holds for~$f$ it also holds for small perturbations of~$f$ (but perhaps with a different~$\eta)$.

For~$p\in\{1,2,\infty\}$, let~$\mu_p$ denote the matrix measure induced by the~$L_p$
vector norm~$|\cdot|_p$. An important advantage of contraction theory is that there exist easy to verify sufficient conditions for contraction in terms of matrix measures. 
For our purposes, it is useful to provide similar conditions for 2-contraction. 
These can be easily derived using the following result.  
\begin{Proposition} (see, e.g.~\cite{muldo1990})
Let~$A\in\R^{n\times n}$. Then
\begin{align*}
 \mu_\infty(A^{[2]})
 & = 
 \max_{1\leq i<j\leq n} \{
 a_{ii}+a_{jj} +\sum_{k\not = i,j } 
 |a_{ik}|+|a_{jk}|  
 \},\\
  \mu_1 (A^{[2]})& = 
 \max_{1\leq i<j\leq n} \{
 a_{ii}+a_{jj} +\sum_{k\not = i,j } 
 |a_{ki}|+|a_{kj}|  
 \},
 \end{align*}
and
$ 
  \mu_2(A^{[2]}) = \lambda_1+\lambda_2,
$
where~$\lambda_1\geq\lambda_2\geq\dots\geq \lambda_n$
are the eigenvalues of~$(A+A^T)/2$.
\end{Proposition}

We say that a dynamical systems has a \emph{non-oscillatory behaviour}~(NOB)
if it has no non-trivial  periodic solutions. In other words, the only possible periodic solutions are equilibrium points. For example, a time-invariant contracting system is~NOB. The same is true for time-invariant~$2$-contracting systems~\cite{muldo1990,li1995}. To illustrate this,
consider the LTI~$\dot x=A x$. If~$A $ is 2-contracting then in particular~$A^{[2]}$ is Hurwitz. Since the eigenvalues of~$A^{[2]}$ are~$\lambda_i+\lambda_j$, $1\leq i<j\leq n$, this implies that~$A$ has no purely imaginary eigenvalues, and thus the LTI is NOB.

Note that the NOB of 2-contracting systems only holds for time-invariant systems. For example, consider the time-varying system:
\begin{align}\label{eq:sint}
    \dot   x_1&=\sin(t),\nonumber\\
    \dot x_2&=-x_2.
\end{align}
The Jacobian of this system is~$J(t,x)=\begin{bmatrix}
0&0\\
0&-1\end{bmatrix}$
and since~$J^{[2]}(t,x)=\trace(J(t,x))\equiv -1$, the system 
is~$2$-contracting.  However, it admits a non-trivial periodic solution, namely,
\[
\begin{bmatrix} x_1(t)\\x_2(t) \end{bmatrix}=
\begin{bmatrix}
1+x_1(0)-\cos(t)\\0
\end{bmatrix},
\] so it is not~NOB. 
The dynamics of~\eqref{eq:sint} contracts   2D polygons  to lines, yet since the system is time-varying, it has a periodic solution along a 1D line. 

Establishing NOB of a dynamical system is important for several reasons. First, certain systems admit a strong Poincar\'e-Bendixson property: any omega limit set that does not include an equilibrium   is a periodic solution.
This holds for example for systems that are monotone with respect to a cone of rank~2~\cite{sanchez2009cones} and in particular for  $3$-dimensional competitive systems~\cite{smith_3d_comp} and for   $2$-cooperative systems~\cite{Eyal_k_posi}.
 If such a system is also NOB then every  omega limit set must contain an equilibrium, and local stability analysis near each equilibrium can often lead to a global picture of the dynamics. This idea has been used
 to provide  a global analysis of many models in epidemiology, see e.g.~\cite{SEIR_LI_MULD1995}. These models are not $1$-contracting, as they typically include two equilibrium points corresponding to the disease-free and the endemic steady states. 
Second,  NOB can sometimes be combined with the closing lemma~\cite{imp_c1} to show that every~$\alpha$ or~$\omega$ limit set of the dynamics consists entirely of equilibria~\cite{li1995}. 

Here, we analyze the serial interconnections of~$k$-contracting systems, with~$k\in\{1,2\}$,
and provide sufficient conditions guaranteeing  that the overall system is NOB or, moreover, 
that every bounded solution converges to an equilibrium.

\begin{figure}
    \centering
    \begin{subfigure}[b]{\columnwidth}
        \centering
        \scalebox{0.9}{
        \begin{tikzpicture}[
            block/.style = {draw, rectangle, thick, minimum height=2em, minimum width=3em}]
            
            \node[block, minimum height=1.25cm, minimum width=2.4cm, label=below:2-contracting] (2contract) {$\begin{aligned}\dot{x}^1 &= f^1(x^1) \\ y &= h(x^1)\end{aligned}$} ;
            \node[block, anchor=north, minimum height=1.25cm, minimum width=2.4cm, label=below:1-contracting] (1contract) [right=1.5cm of 2contract] {$\dot{x}^2 = f^2(x^2, u)$} ;
            \draw[->, thick] (2contract.east) node[above right]{$y$} -- (1contract.west) node[above left]{\raisebox{1.9pt}{$u$}};
        \end{tikzpicture}
        }
    \end{subfigure}%
    \hfill
    \vspace{1em}
    \begin{subfigure}[b]{\columnwidth}
        \centering
        \scalebox{0.9}{
        \begin{tikzpicture}[
            block/.style = {draw, rectangle, thick, minimum height=2em, minimum width=3em}]
            
            \node[block, minimum height=1.25cm, minimum width=2.4cm, label=below:1-contracting] (1contract) {$\begin{aligned}\dot{x}^1 &= f^1(x^1) \\ y &= h(x^1)\end{aligned}$} ;
            \node[block, anchor=north, minimum height=1.25cm, minimum width=2.4cm, label=below:2-contracting] (2contract) [right=1.5cm of 1contract] {$\dot{x}^2 = f^2(x^2, u)$} ;
            \draw[->, thick] (1contract.east) node[above right]{$y$} -- (2contract.west) node[above left]{\raisebox{1.9pt}{$u$}};
        \end{tikzpicture}
        }
    \end{subfigure}%
    \hfill
    \vspace{1em}
    \begin{subfigure}[b]{\columnwidth}
        \centering
        \scalebox{0.9}{
        \begin{tikzpicture}[
            block/.style = {draw, rectangle, thick, minimum height=2em, minimum width=3em}]
            
            \node[block, minimum height=1.25cm, minimum width=2.4cm, label=below:2-contracting] (2contract_1) {$\begin{aligned}\dot{x}^1 &= f^1(x^1) \\ y &= h(x^1)\end{aligned}$} ;
            \node[block, anchor=north, minimum height=1.25cm, minimum width=2.4cm, label=below:2-contracting] (2contract_2) [right=1.5cm of 1contract] {$\dot{x}^2 = f^2(x^2, u)$} ;
            \draw[->, thick] (2contract_1.east) node[above right]{$y$} -- (2contract_2.west) node[above left]{\raisebox{1.9pt}{$u$}};
        \end{tikzpicture}
        }
    \end{subfigure}
    \caption{Serial connections of $k$-contracting systems studied in this note. }
    \label{fig:configurations}
\end{figure}
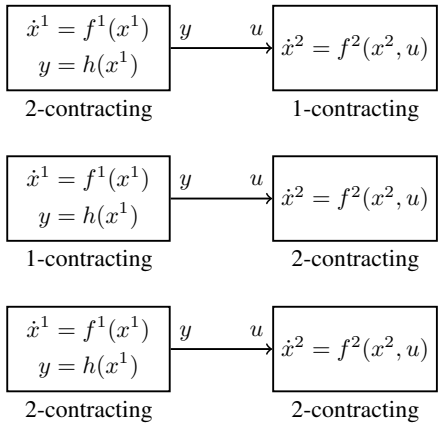
 
\section{Main Results} \label{sec:main}
We begin by studying a serial connection of two sub-systems in the configurations shown in Fig.~\ref{fig:configurations}. We then turn to consider a more general question, namely, when can~$\dot x=f(t,x)$ be decomposed as the serial connection of two systems? We provide a sufficient condition stated as a uniform  ``reducibility condition'' on the Jacobian of~$f$. We then combine these ideas to provide sufficient conditions for well-ordered behaviour of the dynamical system. 
\subsection{Serial connections of two $k$-contracting systems, with~$k\in\{1,2\}$}\label{subsec:ser_con}
Consider the serial  interconnection of two  nonlinear sub-systems.
The first is the time-invariant  sub-system
\begin{align}\label{eq:twosysx}
\dot x^1&=f^1(x^1),\nonumber\\
 y&=h(x^1),
\end{align}
with state~$x^1$ and output~$y$. We assume that the trajectories of this sub-system evolve on a compact and convex set~$\Omega^1$, and that the output map~$h$ is continuous. The second  sub-system is
\be\label{eq:twosysz}
\dot x^2 =f^2(x^2,u),
\ee
with state~$x^2$ and input~$u$. We assume that for any admissible control
the trajectories of this sub-system evolve on a 
compact and convex set~$\Omega^2$.

The  interconnection of the two sub-systems is via~$u=y$ (we assume that~$u,y$ have the same dimension and the same range of admissible values). 
The overall system is thus 
\begin{align}\label{eq:oversys}
\dot x^1&=f^1(x^1),\nonumber\\
\dot x^2&=f^2(x^2,h(x^1)).
\end{align}

The next  two results guarantee the well-ordered asymptotic behaviour 
of the serial connection~\eqref{eq:oversys}. The first result guarantees convergence to an equilibrium (that is not necessarily unique). 

\begin{Proposition}\label{prop:2contract_to_1contract_converge}
Suppose that~\eqref{eq:twosysx} is 2-contracting and that~\eqref{eq:twosysz}
is $1$-contracting.
Then any   solution of~\eqref{eq:oversys}
converges to an equilibrium. 
\end{Proposition}

\begin{Proof}
Fix~$x^i(0)\in\Omega^i$, $i=1,2$. Since~\eqref{eq:twosysx} is 2-contracting, time-invariant, and its trajectories evolve on a compact and convex set, every solution converges to an equilibrium.  Thus,  the limit~$a :=\lim_{t\to\infty} h(x^1(t,x^1(0) ))$  exists. Let~$u_a$ denote
the constant  control~$u(t)\equiv a$.
Since~\eqref{eq:twosysz} is contracting and its trajectories evolve on a compact and convex set, every solution of the system
$\dot x^2 =f^2(x^2,u_a)$  converges to a GAS  equilibrium~$e(a)$.  
This implies that the~$x^2$ system satisfies the converging-input converging-state~(CICS) 
property (see, e.g.,~\cite{cics2006}), so 
$
\lim_{t\to \infty}  x^2( t, x^2(0),h(x^1 (t,x^1(0))))
$
exists, and this completes the proof. 
\end{Proof}

The next result guarantees 
the  non-existence of non-trivial periodic  solutions in the serial connection~\eqref{eq:oversys}.

\begin{Proposition}\label{prop:connection}
Suppose
that~\eqref{eq:twosysx}   is~NOB, and that for any constant input~$u$ system~\eqref{eq:twosysz}
 is~NOB.
  Then~\eqref{eq:oversys}
  is~NOB.
\end{Proposition}

\begin{Proof}
Suppose that~\eqref{eq:oversys} admits a~$T$-periodic solution~$\gamma(t)=\begin{bmatrix} \gamma^1(t)\\\gamma^2(t) \end{bmatrix}$, with~$\gamma(T)=\gamma(0)$. 
Since~\eqref{eq:twosysx}    has no non-trivial
  periodic solutions,~$\gamma^1(t)\equiv  \gamma^1(0)$, so~$h(\gamma^1(t))\equiv h(\gamma^1(0))$.
Thus, along~$\gamma(t)$ the~$x^2$ system has a constant input~$u(t)\equiv h(\gamma^1(0))$ and therefore~$\gamma^2(t)\equiv \gamma^2(0)$. We conclude that~$\gamma(t)\equiv \gamma(0)$, i.e.~$\gamma(t)$ is a trivial periodic solution.
\end{Proof}

It is straightforward to provide conditions guaranteeing that the sub-systems~\eqref{eq:twosysx} and~\eqref{eq:twosysz}  satisfy the requirements  in Prop.~\ref{prop:connection}.
For example, this will be the case if the~$x^1$ system is $2$-contracting,
and the~$x^2$ system is~$2$-contracting for any constant input. 
Note that the~$x^2$ system in~\eqref{eq:oversys}
has a \emph{time-varying} vector field, as~$x^1$ depends  on time.  Still, we can rule out nontrivial  periodic solutions~$\gamma$, because we assume that  along  any such   solution  the component~$\gamma^1(t)$ is constant. 

\begin{figure}[t]
 \begin{center}
   \includegraphics[scale=0.45]{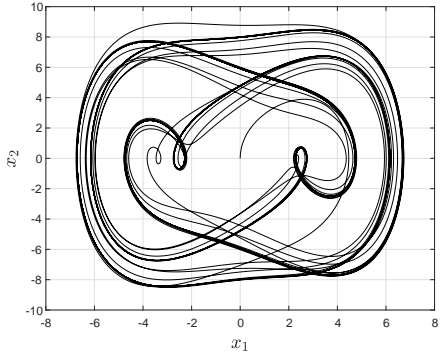}
\caption{Trajectory of the Duffing oscillator with a sinusoidal forcing. }\label{fig:duff_ca}
\end{center}
\end{figure}

Note, however, that even if a time-varying system is NOB for any constant input it may still  display a complicated behaviour for a non-constant input. The next example illustrates  this.
\begin{Example}
Consider the   forced Duffing oscillator
\be\label{eq:ndf}
\ddot x+ \delta \dot x+\beta x^3-\alpha x=\gamma \cos(\omega t).
\ee
with~$ \alpha, \beta,  \gamma,\omega  \geq  0 $ and~$ \delta  > 0 $.
Here the term~$\delta \dot x$ represents a damping term,
and~$\beta x^3$ is a nonlinear restoring force. 
Write~\eqref{eq:ndf} as 
\begin{align}\label{eq:duff}
    \dot x_1&=  x_2  ,\nonumber  \\
    \dot x_2&= -\beta x_1^3 + \alpha x_1 -\delta x_2  +\gamma \cos(\omega t).
\end{align}
The trace of the Jacobian of~\eqref{eq:duff} is~$-\delta<0$, so this  
system is~$2$-contracting and thus~NOB for any \emph{constant} forcing. Fig.~\ref{fig:duff_ca}
depicts the trajectory of~\eqref{eq:duff} for the parameters~$\alpha=0$, $\beta=0.1$,
$\delta=0.1$, $\omega=1$, $\gamma=5$ and
the initial condition~$x_1(0)=x_2(0)=0$. 
It may be seen  that
the trajectory converges to a strange attractor.
\end{Example}

Props.~\ref{prop:2contract_to_1contract_converge} and~\ref{prop:connection} can also be applied to a hierarchical combination of more than two sub-systems. For Prop.~\ref{prop:2contract_to_1contract_converge}, note that the serial connection of any number of 1-contracting sub-systems yields a 1-contracting system, so Prop.~\ref{prop:2contract_to_1contract_converge} may be applied to system where a 2-contracting sub-system feeds a serial connection of multiple 1-contracting sub-systems. For Prop.~\ref{prop:connection}, it is clear from the proof that any number of sub-systems may be used, as long as each sub-system is NOB for any constant input.

Often we are given a dynamical system in the form
\begin{equation} \label{eq:tv_system}
    \dot{x}(t) = f(t,x),
\end{equation}
and not necessarily the interconnection of two sub-systems. We may  still be able to apply Props.~\ref{prop:2contract_to_1contract_converge} and~\ref{prop:connection} above if we can first decompose the given system as the serial connection of two sub-systems. 
\subsection{Decomposing a given system as a serial connection of two sub-systems}
An  interesting and nontrivial  problem is, given an~$n$-dimensional  dynamical system
in the form~\eqref{eq:tv_system},
  can the system be decomposed  as the serial connection of two sub-systems? 
We address this question using a decomposition of~$\R^n$ into two orthogonal subspaces, and a ``uniform reducability'' 
condition on the Jacobian of~$f$.

We assume throughout 
that~$f $ is~$C^1$,   that the solutions of~\eqref{eq:tv_system} evolve on a   convex state-space~$\Omega \subseteq \mathbb{R}^n$, and that for any initial condition~${a \in \Omega}$, and time~$t_0\geq 0$ a unique solution~$x(t,t_0,a)$ exists and satisfies~$x(t,t_0,a) \in \Omega$ for all~$t \geq t_0$.  Let~$J(t,x) := \frac{\partial}{\partial x} f(t,x)$.

Consider an orthogonal decomposition of $\mathbb{R}^n$ into two linear subspaces~$\mathcal{U}$ and~$\mathcal{V}$ of dimensions~$p$ and~$q$, respectively, with~$p,q\geq 1$ and~$p + q = n$. 
Similar decompositions have been used in the context of contraction to subspaces~\cite{Pham2007} or manifolds~\cite{manchester2017control}. The subspaces~$\mathcal{U}$ and~$\mathcal{V}$ are spanned by the columns of~$U \in \mathbb{R}^{n \times p}$ and~$V \in \mathbb{R}^{n \times q}$, respectively, which in turn are chosen such that
\begin{equation}\label{eq:decvu}
    U^T U = I_p, \; V^T V = I_q, \; U^T V = 0_{p \times q},  \; UU^T + VV^T = I_n.
\end{equation}

The next result provides a sufficient condition guaranteeing that~\eqref{eq:tv_system} can be decomposed as the serial connection of two sub-systems.  
\begin{Proposition}\label{prop:decompose}
    Assume that any
 one of the following four equivalent conditions holds: 
 \begin{enumerate}[(a)]
     \item \label{cond:vtf}
    $
        V^T f(t,x) = V^T f(t,VV^Tx)$
          for all $t\in\R_+$, $x\in\Omega $;
\item \label{cond:vtju}
$
V^T J(t,x) U = 0$  for all $t \in \R_+$, $ x \in \Omega$;
\item\label{cond:juu}  
$J(t,x)\mathcal U\subseteq \mathcal U $ for all $ t \in \R_+$, $ x \in \Omega$;
\item\label{cond:jvv}  
$J^T(t,x)\mathcal V\subseteq \mathcal V $ for all $ t \in \R_+$, $ x \in \Omega$.
\end{enumerate}
Let~$y^1:=V^T x$ and~$y^2:=U^Tx$. Then 
   \begin{align}\label{eq:ser2}
   \dot y^1&=V^Tf(t,Vy^1),\nonumber\\
   \dot y^2&= U^T f(t,U y^2+Vy^1).
   \end{align}
\end{Proposition}

Thus, any one of the four equivalent conditions guarantees that~\eqref{eq:tv_system} can be decomposed  as the serial interconnection of the two sub-systems in~\eqref{eq:ser2}, where the output~$Vy^1$ of the~$q$-dimensional~$y^1$ system is fed into the~$p$-dimensional~$y^2$ system. Note that 
 conditions~\eqref{cond:juu} and~\eqref{cond:jvv} are a form of ``uniform
reducibility'' assumption  on the Jacobian of the vector field or its transpose. 

A typical case where the reducibility condition holds is when the dynamics is time-invariant and admits  a first integral~$H(x)=c^T x$, with~$c\in \R^n \setminus \{0\}$. Then along solutions of the system, we have
$
0=\frac{d}{dt} H(x(t)) = c^T f(x(t)),
$
so condition~\eqref{cond:vtf} holds for~$V={c}/{|c|}$.

Note also that
$
Vy^1=VV^T x=x-UU^Tx,
$
that is, the difference between~$x$ and the   (Euclidean norm) projection   of~$x$ on~$\cal U$. 

\begin{Proof}
We first show that the four conditions are equivalent. 
Suppose that condition~\eqref{cond:vtf} holds. Differentiating this condition with respect to~$x$ gives
$
V^T J(t,x)=V^T J(t,VV^T x)V V^T, 
$
and multiplying by~$U$ on the right  gives~\eqref{cond:vtju}.
To prove the converse implication,
 assume that~\eqref{cond:vtju} holds. Then
\begin{equation}
\begin{split}
    0 &= \int_0^1 V^T J(t,VV^Tx + rUU^Tx) UU^Tx \mathrm{d}r \\
      &= \int_0^1 V^T \frac{\partial}{\partial r} f(t,VV^Tx + rUU^Tx) \mathrm{d}r \\
      &= V^T f(t,x) - V^T f(t,VV^T x),
\end{split}
\end{equation}
where in the last equation we used~\eqref{eq:decvu}. Thus,~\eqref{cond:vtju} implies~\eqref{cond:vtf}.
The equivalence of~\eqref{cond:vtju},~\eqref{cond:juu}, and~\eqref{cond:jvv} follows from~\eqref{eq:decvu}.
Now suppose that condition~\eqref{cond:vtf} holds.
   Let~$y(t): = \begin{bmatrix}
    y^1(t)\\y^2(t)
 \end{bmatrix} $,
    with~$y^1:=V^T x$ and~$y^ 
    2 : =U^T x$. 
    Then  
    \begin{align}\label{eq:y1dyn}
        \dot y^1&= 
        V^T f(t,x) 
 \nonumber\\&=
        V^T f(t,VV^Tx) 
        \nonumber \\&=
        V^T f(t,Vy^1), 
    \end{align}
and 
 \begin{align*}
 \dot y^2&=U^T f(t,x)\\&=U^T f (t,  UU^T x + VV^T  x)\\
 &=U^T f (t,U y^2+ V y^1),
 \end{align*}
  and this completes the proof.
\end{Proof}

\begin{Example}
Consider  the nonlinear system
\[
\dot x=f( L  x )  ,
\]
where~$ L$ is the Laplacian
of a weighted digraph. 
For example, if~$f(y)=-y $ we get the linear consensus protocol, whereas if~$f(y)=\begin{bmatrix} -\tanh(y_1)&\dots&-\tanh(y_n) \end{bmatrix}$   we get a  form  of a ``bounded derivatives'' consensus protocol.

Since~$L 1_n=0$, we can take~$U:=\frac{1}{\sqrt{n}}1_n$, and let~$V$ be as in~\eqref{eq:decvu}. Then for~$y^1:=V^T x$, $y^2:=U^Tx$, we have
\[
\dot y^1 =  V^T  f( Lx)=
V^T f( L(U U^T +V V^T) x ) 
=V^T f( L V y^1),
\]
and
\[
\dot y^2= (1/\sqrt{n})1_n^T f( L x )    =(1/\sqrt{n})
 1_n^T f(L  V y^1).
\]
The~$(n-1)$-dimensional~$y^1$ system describes the  dynamics on the subspace  orthogonal to the ``consensus subspace''~${\cal U} = \spanop(1_n)  $.
The dynamics of the~$1$-dimensional~$y^2$  system depends on
$
    Vy^1 ,
$
that is, the difference between~$x$ and its (Euclidean norm) projection on~$\cal U$.
\end{Example}
 
\begin{Remark}
Prop.~\ref{prop:decompose} implies the well-known result that any LTI system~$\dot x=A x$, $A \in \R^{n \times n}$, with $n \geq 3$, may be decomposed into a serial interconnection of two sub-systems. If $A$ has a real eigenvalue~$\lambda$, with a  corresponding real eigenvector~$u$,  then~$\mathcal{U}$ may be chosen as the subspace spanned by~$u $. Otherwise, $A$ has a pair of complex conjugate eigenvalues~$\alpha \pm j \beta$ and corresponding eigenvectors $u^1 \pm j u^2$, where $\alpha \in \R$, $\beta \in \R\setminus\{0\}$, $u^1  \in \R^n$, $ u^2 \in \R^n\setminus\{0\}$, and~$j=\sqrt{-1}$. Let~$\mathcal{U}:= \spanop(u^1, u^2)$. Then 
for any $z \in \mathcal{U}$,
\[
A z = A (r u^1 + s u^2) = (r \alpha + s \beta)u^1 + (s\alpha - r\beta)u^2 , \]
so~$A$ maps~$\mathcal{U}$ to~$\mathcal{U}$.
\end{Remark}

\begin{Remark}
The conditions in Prop.~\ref{prop:decompose} are related to the existence of invariant sets of the dynamics~\eqref{eq:tv_system} or of its corresponding \emph{variational equation}. Indeed, pick two initial conditions~$a,b \in\Omega$, an initial time~$t_0:=0$, 
and let~$z(t):=x(t,a)-x(t,b)$.
Then
\begin{align*}
\dot z(t) &= f(t,x(t,a)) -f(t,x(t,b)) \\
&=\int_0^1 \frac{\partial}{\partial r}
f( t, r x(t,a)  +(1-r) x(t,b)) \diff r\\
&=\int_0^1 J(t,  r x(t,a)  +(1-r) x(t,b))\diff r\, z(t).
\end{align*}
Combining this variational equation with condition~\eqref{cond:juu} implies that if~$a-b\in \cal U$ then~$x(t,a)-x(t,b) \in \cal U$ for all~$t\geq 0$.

If we   assume, in addition,  that~$V^Tf(t,0)=0$ for all $t \ge 0$ (this holds, for example,  if~$0$ is an equilibrium  of~\eqref{eq:tv_system} for 
all~$t \ge 0$)
then condition~\eqref{cond:vtf} gives
\[
    V^T f(t,U U^T x) = V^T f(t, VV^T  U U^T x       ) = V^Tf(t,0) \equiv 0  ,
\]
and this implies that 
$
    f(t,\mathcal{U}) \subseteq  \mathcal{U},
$
i.e.~$\mathcal{U}$ is an invariant set of the dynamics.
\end{Remark}
 
The next result demonstrates an application of 
Prop.~\ref{prop:decompose} to a system with a feedback form.
\begin{Corollary}
\label{cor:second_order_consensus}
Consider an~$(n+m)$-dimensional   system:
\begin{equation}\label{eq:yoash_example}
\begin{aligned}
    \dot{x} &= g(z), \\
    \dot{z} &= h(M x,z),
\end{aligned}
\end{equation}
where~$x\in\R^n$,
$ z \in\R^m$, and~$M \in  \R^{k\times n}$. Suppose that
there exist~$U\in\R^{n\times p}$ and~$V\in\R^{n\times (n-p)}$
as in~\eqref{eq:decvu} such that
\be\label{eq:muzero}
MU=0.
\ee
Define~$y^1:=V^Tx$,
$y^2:=z $, and~$y^3:=U^T x$. 
Then
\begin{align}\label{eq:splsys}
\dot y^1&=V^T g(y^2) \nonumber,\\
\dot y^2&=h(M V y^1,y^2) \nonumber,\\
\dot y^3&=U^T g(y^2) .
\end{align}
\end{Corollary}
Note that this implies   a decomposition into an~$(n+m-p)$-dimensional  sub-system with state~$\begin{bmatrix} y^1\\y^2\end{bmatrix}$, whose output~$g(y^2)$ is fed into the~$p$-dimensional~$y^3$ sub-system.
Note also that if~$\rank(M)<n$ then we can always find~$U,V$ satisfying  condition~\eqref{eq:muzero}.

\begin{Proof}
Define matrices~$\tilde U \in \R^{(n+m) \times p}$  and~$\tilde V \in \R^{(n+m)\times (n+m-p) }$
by
$
\tilde U := \begin{bmatrix} U\\ 0\end{bmatrix}$,
$ 
\tilde V := \begin{bmatrix} V&0\\ 0&I\end{bmatrix}.
$
The Jacobian of~\eqref{eq:yoash_example} is
\begin{equation}
    J(x,z) = \begin{bmatrix}
        0 &  \frac{\partial}{\partial z} g(z) \\
     \frac{\partial}{\partial x}h(Mx,z)M & \frac{\partial}{\partial z}h(Mx,z)
    \end{bmatrix},
\end{equation}
so~$J(x,z) \begin{bmatrix} 
U\\ 0\end{bmatrix} = 0$ for all~$x,z$. Thus, $J$ maps the subspace~${ \cal \tilde U} $ corresponding to~$\tilde U$ 
to itself.
By 
Prop.~\ref{prop:decompose}, 
the system can be decomposed as the serial connection of two sub-systems by defining~$\tilde y^1:=\tilde V^T \begin{bmatrix} x\\ z\end{bmatrix}$ and~$\tilde y^2:=\tilde U^T \begin{bmatrix} x\\ z\end{bmatrix}$, that is,~$\tilde y^1 = \begin{bmatrix} V^T x\\ z\end{bmatrix}$ and~$\tilde y^2:=  U^T  x $. Computing the dynamics of $y^1$ and $y^2$ yields~\eqref{eq:splsys}.
\end{Proof}

\begin{Example}
 \label{exa:2nd_2}
Consider the second-order consensus system~\cite{2nd_order_cons}:
\begin{align}\label{eq:2od_2}
\dot x&=v   ,\nonumber \\
\dot v&=u,
\end{align}
with
\[
u= f( \beta L x +\alpha L v).
\]
Here~$ x\in\R^n$, $x_i$ describes the (scalar) 
location of     agent~$i$,~$v_i$ is the velocity of agent~$i$,~$L$ is the Laplacian
of a weighted digraph, and $\alpha,\beta>0$. The goal is to drive both the~$x_i$s and the~$v_i$s to consensus, that is,~$x(t)\to c_1 1_n$, $v(t)\to c_2 1_n$ for some constants~$c_1,c_2$. 
The nonlinear 
function  can be for example~$f(y)=-y$ or~$f(y)=-\tanh(y)$. 

This system has the form  described in Corollary~\ref{cor:second_order_consensus}
with~$z=v$, $m=n$, $g(z)=z$, $M=\beta L$, 
and~$h(Mx,z)=f(Mx+\alpha L z)$. We can thus take
$U=\frac{1}{\sqrt{n}}1_n$,
and a~$V\in\R^{n\times(n-1)}$ such that~\eqref{eq:decvu} holds. Then for~$y^1:=V^T x$,
$y^2:=v$, and~$y^3:=\frac{1}{\sqrt{n}}1_n^Tx$, we get
\begin{align*}
\dot y^1&=V^T y^2,\\ 
\dot y^2&=f(\beta L V y^1+\alpha L v),\\ 
\dot y^3&= U^T y^2 . 
\end{align*}
The first two equations form a~$(2n-1)$ dimensional system, and the third describes the dynamics of~$1_n^T x$. 
\end{Example}

In the next subsection, we combine  Props.~\ref{prop:2contract_to_1contract_converge} and~\ref{prop:connection}  with the reducibility condition to guarantee the well-ordered behaviour of
the time-invariant system 
\be\label{eq:tif}
\dot x=f(x).
\ee
We assume throughout  that~$f$ is~$C^1$ and that the trajectories evolve on a convex set~$\Omega \subseteq \R^n$.
Let~$J(x):=\frac{\partial}{\partial x} f(x)$.

 \subsection{Conditions for well-ordered behaviour of~$\dot x=f(x)$}
The next result provides a sufficient condition for NOB  that is based on 2-contraction on a certain subspace. As shown in Example~\ref{exa:lti} below, this is weaker than requiring 2-contraction on the entire state-space. 

\begin{Proposition}\label{prop:V_2contraction}
    Suppose that the decomposition condition in Prop.~\ref{prop:decompose} holds for~\eqref{eq:tif} with~$\mathcal{U}$   a  one-dimensional subspace, i.e. $p = 1$, and that  there exists a matrix measure $\mu$ such that
    \begin{equation} \label{eq:V_2contract}
        \mu\left((V^T)^{(2)} J^{[2]}(x) V^{(2)} \right ) \le -\eta < 0, \text{ for all } x\in\Omega.
    \end{equation}
   Then \eqref{eq:tif} is~NOB.
\end{Proposition}

\begin{Proof}
We know that 
for~$y^1:=V^T x$ and~$y^ 
    2 : =U^T x$, we have~$
        \dot y^1 = 
         V^T f(Vy^1), 
$ and
 $\dot y^2  =U^T f (U y^2+ V y^1)$. The Jacobian of the~$y^1$ system is
 $V^T J(V y^1) V$, and combining this with~\eqref{eq:ctf} and~\eqref{eq:V_2contract} implies that the~$y^1$ system is 2-contracting and thus~NOB. 
The~$y^2$ system
 is   one-dimensional   and thus~NOB for any constant input.
 Applying  Prop.~\ref{prop:connection}   completes the proof. 
\end{Proof}

Note that the existence of  a one-dimensional invariant subspace 
is quite common in various systems, e.g. in models for  synchronization, where the synchronized state (i.e.,~$ x_1=\dots=x_n$) is invariant, see also the examples in Section~\ref{sec:app}. 

It is instructive to demonstrate Prop.~\ref{prop:V_2contraction}
in the case of an  LTI  system.
\begin{Example}\label{exa:lti}
Consider the LTI system
\begin{equation} \label{eq:lti_system}
    \dot{x}(t) = A x(t).
\end{equation}
In this case, the decomposition  condition is~$A\, {\mathcal{U}}\subseteq \mathcal{U}$.
Let~$T$ be the~$n\times n$ matrix
\be  \label{eq:tit}
    T := \begin{bmatrix}
        V^T \\
        U^T \\
    \end{bmatrix}.
\ee
Then~$   T^{-1} = \begin{bmatrix}
       V&U
    \end{bmatrix}$, and
\begin{equation}\label{eq:tat}
    T A T^{-1} = \begin{bmatrix}
        U^T A U & U^T A V \\
        0       & V^T A V
    \end{bmatrix}.
\end{equation}
Thus,~$TAT^{-1}$ is reducible. The spectrum of $A$ is the union of
the real scalar 
 $U^T A U$ and the spectrum of~$V^T A V$.   Since~$V^T A V$ is 2-contracting, it  has  no pure imaginary eigenvalues, so   $A$ has no pure imaginary eigenvalues. Thus,~\eqref{eq:lti_system}
 has no  non-trivial periodic trajectories. 

It is important to note that the eigenvalue~$U^TAU$ of~$TAT^{-1}$, and thus of~$A$,
can be arbitrarily large, 
so~$A$ is not necessarily $2$-contracting in the entire state-space. 
\end{Example}

Roughly speaking, 
Prop.~\ref{prop:V_2contraction}
requires that~$\mathcal U$ is one-dimensional and that the system is 2-contracting  on~$\mathcal V$, and proves that such a configuration is NOB. By requiring that the system is instead 1-contracting on~$\mathcal U$ (no longer necessarily one-dimensional), we now derive a stronger result, namely that every bounded trajectory of the overall system converges to an equilibrium point. Example~\ref{exa:lti2exa} below shows that these conditions do not imply that the system is 2-contracting on the entire state-space.

\begin{Proposition}\label{prop:thm2}
    Suppose that the decomposition condition in Prop.~\ref{prop:decompose} holds, and that
    there exists a matrix measure $\mu$ such that
    \begin{equation} \label{eq:V_2contracting}
        \mu((V^T)^{(2)} J^{[2]}(x) V^{(2)}) \le -\eta < 0
    \end{equation}
    and
    \begin{equation}\label{eq:U_1contracting}
        \mu(U^T J(x) U) \le -\beta < 0
    \end{equation}
    for all $x \in \Omega$. Then every bounded trajectory of \eqref{eq:tif} converges to an equilibrium point.
\end{Proposition}

\begin{Proof}
For~$y^1:=V^T x$ and~$y^ 
    2 : =U^T x$, we have~$
        \dot y^1 = 
         V^T f(Vy^1), 
$ and
 $\dot y^2  =U^T f (U y^2+ V y^1)$. The Jacobian of the~$y^1$ system is
 $V^T J(V y^1) V$, and combining this with~\eqref{eq:ctf} and~\eqref{eq:V_2contracting} implies that the~$y^1$ system is 2-contracting.
 The Jacobian of the~$y^2$ system is $U^T J(Vy^1 + Uy^2) U$, so~\eqref{eq:U_1contracting} implies that the~$y^2$ system is 1-contracting. Applying Prop.~\ref{prop:2contract_to_1contract_converge} completes the proof. 
\end{Proof}

Again, it is instructive to demonstrate Prop.~\ref{prop:thm2}
in the case of an LTI system.
\begin{Example}\label{exa:lti2exa}
Consider the LTI system~\eqref{eq:lti_system}. The decomposition condition implies that~$V^T A U=0 $, so~\eqref{eq:tat} holds. 
Eq.~\eqref{eq:U_1contracting} implies that all the eigenvalues of~$U^TAU$ have a negative real part.
Eq.~\eqref{eq:V_2contracting} implies that~$V^TAV$ has no pure imaginary eigenvalue. We conclude that  the spectrum
of~$A$ has no pure imaginary eigenvalues. Thus, any bounded trajectory of the LTI converges to an equilibrium point.

Note that the conditions do not imply that the overall system is 2-contracting on the entire state-space. Consider for example the LTI system~\eqref{eq:lti_system} with
$
    A = \diag(2,-3,-1,-1)
$. Let
$
    U := \begin{bmatrix}
        0 \\
        I_2 
    \end{bmatrix},
    V := \begin{bmatrix}
        I_2 \\
        0 
    \end{bmatrix}.
$
Then, for any monotonic norm $\mu(U^T A U) \le -1$ and $\mu((V^T)^{(2)} A^{[2]} V^{(2)}) = -1$. The decomposition condition also holds. However, the maximal  eigenvalue of $A^{[2]}$ is one, so the system is not 2-contracting on the entire state-space for any norm.
\end{Example}

\section{ Applications}\label{sec:app}
We describe two  simple applications  of the theoretical results. We will make use of the following fact (see, e.g.~\cite{schwarz1970}). 
If~$A\in\R^{3\times 3}$  then  
\be\label{eq:a2comp}
A^{[2]} = \begin{bmatrix} 
a_{11}+a_{22}&a_{23} &-a_{13} \\
a_{32} & a_{11}+a_{33} &a_{12}\\
-a_{31} & a_{21} & a_{22}+a_{33}
\end{bmatrix} .
\ee

Our first application is a 3D system with two agents.
\begin{Corollary}\label{coro:3d}
Consider the system
\begin{align}\label{eq:3D}
\dot x_1&=f(x_1,x_2),\nonumber \\
\dot x_2&=x_3- x_1,\\
\dot x_3&=g(x_3,x_2) .\nonumber
\end{align}
Assume that
\be\label{eq:firass}
\frac {\partial}{\partial x_1} f(x_1,x_2)=\frac {\partial}{\partial x_3} g(x_3,x_2)\leq -\delta< 0 \text{ for all }x,
\ee
and that the trajectories evolve on a convex and compact set. Then every trajectory of~\eqref{eq:3D} converges to an equilibrium point.
\end{Corollary}

Here~$x_1$ and~$x_3$ may  represent the state of  two ``agents'', 
and~$x_2$ evolves according to the difference between the agent states. 
A typical example is a system describing the interconnection of two synchronous generators, that interact via an integral of the difference between their frequencies (i.e. the relative phase angle)~\cite{Kundur1994}. In the control theory community, such models are often called network reduced power systems.

\begin{Proof}
The Jacobian of~\eqref{eq:3D} is
$
J =\begin{bmatrix}
\frac {\partial}{\partial x_1} f&\frac {\partial}{\partial x_2} f &0 \\
 -1&0& 1 \\
0&
\frac {\partial}{\partial x_2} g&\frac {\partial}{\partial x_3}  g 
\end{bmatrix}.
$
Consider the orthogonal decomposition generated by the matrices
\[
    U := \frac{1}{\sqrt{2}} \begin{bmatrix}
        1 \\ 0 \\ 1
    \end{bmatrix}, \quad 
    V := \frac{1}{\sqrt{2}}\begin{bmatrix}
        0        & 1  \\
        \sqrt{2} & 0  \\
        0        & -1 \\
    \end{bmatrix}.
\]
Then
$
V^T J U= (1/2) \begin{bmatrix} 
0\\\frac {\partial}{\partial x_1} f-\frac {\partial}{\partial x_3} g
\end{bmatrix} ,
$
and~\eqref{eq:firass} gives~$V^T J U=0$.
Additionally, 
$V^{(2)}=\begin{bmatrix}
-1&0&-1\end{bmatrix}^T $,
and
\[
J ^{[2]}=\begin{bmatrix}
\frac {\partial}{\partial x_1} f & 1 &0 \\
\frac {\partial}{\partial x_2} g & \frac {\partial}{\partial x_1} f+ \frac {\partial}{\partial x_3} g& \frac {\partial}{\partial x_2} f \\
0&-1&
\frac {\partial}{\partial x_3} g
\end{bmatrix},
\]
so~\eqref{eq:firass} yields 
$
(V^T)^{(2)} J^{[2]}  V^{(2)}
=
\frac {\partial}{\partial x_1} f+\frac {\partial}{\partial x_3} g \leq -2\delta < 0$, and the 2D sub-system is 2-contracting. Similarly, $U^T J U \leq -2\delta < 0$ and the scalar sub-system is 1-contracting. Applying Prop.~\ref{prop:thm2} completes the proof of Corollary~\ref{coro:3d}. 
\end{Proof}

Our second application describes a system of three ``synchronizing agents''.
\begin{Corollary}\label{coro:cos}
 Consider the system:
\begin{align}\label{eq:nopos}
    \dot x_1 &= f(x_1-x_2) + f( x_1-x_3)  ,\nonumber\\
    \dot x_2  &=   g(x_2-x_1) + g( x_2-x_3) , \\
    \dot x_3 &=  h(x_3-x_2) + h( x_3-x_1)   \nonumber,
\end{align}
where~$f,g,h:\R\to\R$ are~$C^1$. 
Suppose that the trajectories evolve on a compact and convex set~$\Omega$ and that 
\be\label{eq:zcond}
f'(p)+g'(-p) + f'(q)+h'(-q)+ g'(q-p)+h'(p-q)\leq-\eta < 0  
\ee
 for any~$p,q\in\R$.
Then~\eqref{eq:nopos} is~NOB. 
\end{Corollary}

\begin{Proof}
 Let~$x_{ij}:=x_i-x_j$.
 The Jacobian~$J(x)$ of~\eqref{eq:nopos} is
\[
 \begin{bmatrix} 
f'( x_{12}) +f'(x_{13}) & -f'(x_{12})  &  -f'(x_{13}) \\
-g'( x_{21} )&  g'(x_{21}) + g'(x_{23} )  &  -g'(x_{23}) \\
-h'( x_{31}) &  -h'(x_{32}) & h'(x_{32})  +   h'(x_{31}) 
\end{bmatrix},
\]
and
using~\eqref{eq:a2comp} implies that~ 
\[ 
 J^{[2]}(x)= \left [  \begin{smallmatrix} 
    J_{11}(x)+J_{22}(x)
 & -g'(x_{23}) & f'(x_{13}) \\
-h'(x_{32}) & 
  J_{11}(x)+J_{33}(x) & -f'(x_{12}) \\
h'(x_{31})  & -g'(x_{21}) &  J_{22}(x)+J_{33}(x)
\end{smallmatrix} \right ].
\]
Let~$U:=\frac{1}{\sqrt{3}} \begin{bmatrix}  1\\1\\1 \end{bmatrix}$,
$V:=\frac{ 1}{\sqrt{6}}  \begin{bmatrix} 2&0\\
-1&-\sqrt{3} \\
-1&\sqrt{3} \end{bmatrix}$.
Then~$J(x)U=0$ for all~$x$, so the reducibility  condition  holds. 
A calculation gives~$V^{(2)}=\frac{1}{\sqrt{3}} \begin{bmatrix}-1&1&-1 \end{bmatrix}^T $, and
\begin{align*}
(V^T)^{(2)}  & J^{[2]})(x) V^{(2) }   = f'(x_{12}) +f'(x_{13})\\&+
g'(x_{21}) +g'(x_{23})+h'(x_{31}) +h'(x_{32}).
\end{align*}
Combining this with~\eqref{eq:zcond} implies that all the conditions
in
Prop.~\ref{prop:V_2contraction} hold. This completes the proof of Corollary~\ref{coro:cos}. 
\end{Proof}

\section{Conclusion}
An important topic in systems theory is analyzing an interconnected system based on the properties of the sub-systems and the interconnection network.  
In this context, an important advantage of contracting systems is that various interconnections of such systems yield a contracting system. 

We analyzed the serial interconnection
of~$k$-contracting systems, with~$k\in\{1,2\}$. 
Our results guarantee NOB and, under stronger assumptions, that every bounded solution converges to an equilibrium (that is not necessarily unique). 
To apply these results to a wider set of systems, we also derived a reducibility condition  guaranteeing that a given system can be decomposed  as the serial connection of two systems.

Prop.~\ref{prop:decompose}
provides a sufficient condition for decomposing a system as the serial connection of two sub-systems based on a decomposition of~$\R^n$ into two subspaces. It may be of interest to  
extend this result using  more general decompositions of~$\R^n$. 

Our reducibility condition is restrictive and not robust to small perturbations  in the dynamics. Another topic for further research is to apply our results to a system~$\dot x=f(x)$ that does not satisfy the reducibility condition using the following scheme: (1) approximate~$f$ using a vector field~$g$ that does satisfy the reducibility  condition; (2) analyze the dynamics~$\dot y=g(y)$ using the tools developed here; and (3) use comparison principles for ODEs~\cite{Szarski1965} 
to show that the results for the~$y$-system also hold for the original~$x$-system. These topics are currently under study.


\end{document}